\documentclass[11pt,oneside]{article}
\usepackage{mathrsfs}
\usepackage{amssymb}
\usepackage{amsfonts}
\usepackage{bbding}
\usepackage{pifont}
\usepackage{fancyhdr}
\usepackage{titlesec}
\usepackage{cite}
\usepackage{ifthen}
\usepackage{amssymb,amsthm,amscd}

\titleformat{\section}{\centering\Large\bfseries}{\S\arabic{section}}{1em}{}
\newboolean{first}
\setboolean{first}{true}
\renewcommand{\headrulewidth}{0pt}

\begin{document}

\setlength\abovedisplayskip{2pt}
\setlength\abovedisplayshortskip{0pt}
\setlength\belowdisplayskip{2pt}
\setlength\belowdisplayshortskip{0pt}

\title{\bf \huge Hyperbolic K${\bf\ddot{a}}$hler-Ricci Flow}
\author{Xu Chao\footnote{e-mail address: xuchaomykj@163.com}\\
        Department of Mathematics\\
        Zhejiang University, Hangzhou, China} \maketitle

%\begin{center}
%Abstract
%\end{center}

\begin{abstract}
%\begin{center}
%\begin{minipage}{100mm}
In this paper, the author has considered the hyperbolic
K$\ddot{a}$hler-Ricci flow introduced by Kong and Liu \cite{KL},
that is, the hyperbolic version of the famous K$\ddot{a}$hler-Ricci
flow. The author has explained the derivation of the equation and
calculated the evolutions of various quantities associated to the
equation including the curvatures. Particularly on Calabi-Yau
manifolds, the equation can be simplified to a scalar hyperbolic
Monge-Amp$\grave{e}$re equation which is just the hyperbolic version
of the corresponding one in K$\ddot{a}$hler-Ricci flow.\\\hskip 2mm
{\small }
%\end{minipage}
%\end{center}
\end{abstract}
\newpage

\thispagestyle{fancyplain} \fancyhead{}
%\fancyhead[L]{\textit{paper}\\
%\fancyfoot{} \vskip 10mm

\section{Introduction}

Recently, a new flow on Riemannian manifolds is introduced by Kong
and Liu(\cite{KL}, \cite{K}, \cite{DKL}):
\\$$\frac{\partial^2}{\partial t^2}g=-2Rc.$$\\
This is the hyperbolic version of the famous Ricci flow. They
studied its short-time existence in compact case, derived evolutions
of curvatures which have wave character. A dissipative flow was also
considered in \cite{DKL2}. One remarkable result is \cite{KLX}, in
which the authors proved that on compact surfaces with metric
depending on one space variable, the flow has a global solution if
the initial velocity is large enough. So given an initial metric,
the geometric surgery can be replaced by choosing suitable velocity
to allow the long-time existence. However, the evolutions of Riemann
curvature, Ricci curvature and the scalar curvature are very
complicated which contain the first-order derivative of time that is
hard to control. In this paper, I consider its complex version which
is also introduced by Kong and Liu in \cite{KL}:
\\$$\frac{\partial^2}{\partial t^2}g_{\alpha\bar{\beta}}=-R_{\alpha\bar{\beta}}.$$\\
As we will see on Calabi-Yau manifolds, it can be reduced to a
scalar (complex) hyperbolic Monge-Amp$\grave{e}$re equation which is
much easier to handle. The flow considered here can also be regarded
as the hyperbolic version of K$\ddot{a}$hler-Ricci flow:
\\$$\frac{\partial}{\partial t}g_{\alpha\bar{\beta}}=-R_{\alpha\bar{\beta}}.$$\\

In the past several decades, nonlinear partial differential
equations played important role in differential geometry. For
example, the resolutions of Calabi's conjecture \cite{Y} and
Poincar$\acute{e}$'s conjecture \cite{CZ} \cite{MT} \cite{P1}
\cite{P2} \cite{P3} are due to the methods of geometric analysis and
the elliptic and parabolic type of equations are intensively
studied. However, the hyperbolic equations have been ignored for a
long time. For the elliptic and parabolic equations, we have the
powerful tool the maximum principle which cannot be applied to
hyperbolic ones. This lack of practical tool may explain the reason
why the research of hyperbolic equations on manifolds are not as
active as the other two types. However, Perelman \cite{P1} \cite{P2}
\cite{P3} introduced several powerful tools into the study of flow
like the energy functional, the monotonicity formula, and the
space-time geometry. In the Euclidean space, these methods are
available for hyperbolic equations. So I believe we can refer to
Perelman's methods to study nonlinear hyperbolic equations on
manifolds.

Hyperbolic equations are important in physics especially in general
relativity. The famous Einstein's equation has been studied for a
long time. In fact, as illustrated in \cite{KL}, the hyperbolic
version of Ricci flow is closely related to Einstein's equation. So
I hope, by studying the wave character of manifolds through
hyperbolic equations, we can better understand the basic structure
of the universe. There is already some work on wave equations on
Lorentzian manifolds which catches more interest from Physicists.
Although wave equations on Riemannian manifolds are also considered,
the results are few and have a more analytic style that its
influence on curvature are not considered. In my opinion, analyzing
wave equations on manifolds from a more geometric viewpoint so that
it may give more help to understand the geometry and topology of
manifolds.\\

Let $(\mathcal{M}^{n}, g, J)$ be a complete K$\ddot{a}$hler
manifold, where $g$ is the metric depending on time and $J$ is a
fixed complex structure. The hyperbolic K$\ddot{a}$hler-Ricci flow
is the following evolution equation
\\$$\frac{\partial^2}{\partial t^2}g_{i\bar{j}}=-R_{i\bar{j}}\eqno(1.1)$$\\
for a family of K$\ddot{a}$hler metrics $g_{i\bar{j}}(t)$ on
$\mathcal{M}^n$. A natural and fundamental problem is the short-time
existence and uniqueness theorem of (1.1). We will see as in the
case of K$\ddot{a}$hler-Ricci flow, the evolution equation of
metrics can be simplified to a scalar evolution equation. In the
present paper, I also derive the corresponding
wave equations for various geometric quantities.\\

{\bf Theorem 1.1. }{\it Let $(\mathcal{M}^n, g^0_{i\bar{j}}(z))$ be
a compact Calabi-Yau manifold with $g^0_{i\bar{j}}(z)$ a
K$\ddot{a}$hler metric. Then there exists a constant $h>0$ such that
the initial value problem
\\$$\cases{
    \frac{\partial^2}{\partial
    t^2}g_{i\bar{j}}(z,t)=-R_{i\bar{j}}(z,t)\cr
    g_{i\bar{j}}(z,0)=g^0_{i\bar{j}}(z),\quad \frac{\partial}{\partial
    t}g_{i\bar{j}}(z,0)=g^1_{i\bar{j}}(z),\cr}
    $$\\
has a unique smooth solution $g_{i\bar{j}}$ on
$\mathcal{M}\times[0,h]$ and $g_{i\bar{j}}(z,t)$ remains
K$\ddot{a}$hler for any $t>0$ as long as the solution exists and the
initial velocity must satisfy $[\frac{\partial}{\partial
t}\omega(z,0)]=0$ where $\omega(z,0)$ is the K$\ddot{a}$hler form of
$g_{i\bar{j}}(z,0)$.}\\

Similar to K$\ddot{a}$hler-Ricci flow, I will derive the
corresponding
wave equations for the curvatures.\\

{\bf Theorem 1.2. }{\it Under the hyperbolic K$\ddot{a}$hler-Ricci
flow (1.1), the Riemannian curvature tensor, Ricci curvature and
scalar curvature satisfy the evolution equations under a unitary
coordinate
\\ \begin{eqnarray*}\frac{\partial^2}{\partial t^2}R_{i\bar{j}k\bar{l}}&=&\triangle_R R_{i\bar{j}k\bar{l}}
+R_{i\bar{\alpha}\beta\bar{l}}R_{\alpha\bar{j}k\bar{\beta}}-R_{i\bar{\alpha}k\bar{\beta}}R_{\alpha\bar{j}\beta\bar{l}}
+R_{i\bar{j}\beta\bar{\alpha}}R_{\alpha\bar{\beta}k\bar{l}}\\
&
&-\frac{1}{2}(R_{i\bar{\alpha}}R_{\alpha\bar{j}k\bar{l}}+R_{\alpha\bar{j}}R_{i\bar{\alpha}k\bar{l}}
+R_{k\bar{\alpha}}R_{i\bar{j}\alpha\bar{l}}+R_{\alpha\bar{l}}R_{i\bar{j}k\bar{\alpha}})\\
& &+2g^{p\bar{q}}\nabla_{k}\biggl(\frac{\partial}{\partial
t}g_{i\bar{q}}\biggr)\nabla_{\bar{l}}\biggl(\frac{\partial}{\partial
t}g_{p\bar{j}}\biggr)\end{eqnarray*}
\\ \begin{eqnarray*}\frac{\partial^2}{\partial t^2}R_{i\bar{j}}&=&\triangle_R
R_{i\bar{j}}+R_{i\bar{j}k\bar{l}}R_{l\bar{k}}-R_{i\bar{k}}R_{k\bar{j}}
-2\biggl\langle\frac{\partial}{\partial
t}g_{k\bar{l}},\frac{\partial}{\partial
t}R_{i\bar{j}k\bar{l}}\biggr\rangle\\
& &+2R_{i\bar{j}k\bar{l}}\biggl(\frac{\partial}{\partial
t}g_{n\bar{m}}\biggr)\biggl(\frac{\partial}{\partial
t}g_{r\bar{s}}\biggr)g^{r\bar{m}}g^{n\bar{l}}g^{k\bar{s}}\\
&
&+2g^{p\bar{q}}g^{k\bar{l}}\nabla_{k}\biggl(\frac{\partial}{\partial
t}g_{i\bar{q}}\biggr)\nabla_{\bar{l}}\biggl(\frac{\partial}{\partial
t}g_{p\bar{j}}\biggr)\end{eqnarray*}
\\ \begin{eqnarray*}\frac{\partial^2}{\partial t^2}R&=&\triangle
R+|R_{\alpha\bar{\beta}}|^2
+\triangle\biggl|\frac{\partial}{\partial
t}g_{\alpha\bar{\beta}}\biggr|^2
-2\biggl\langle\frac{\partial}{\partial t}g_{\alpha\bar{\beta}},
\frac{\partial}{\partial
t}R_{\alpha\bar{\beta}}\biggr\rangle\\
& &+2R_{k\bar{l}}\biggl(\frac{\partial}{\partial
t}g_{n\bar{m}}\biggr)\biggl(\frac{\partial}{\partial
t}g_{r\bar{s}}\biggr)g^{r\bar{m}}g^{n\bar{l}}g^{k\bar{s}}.\end{eqnarray*}\\
Here
$\triangle_R=\frac{1}{2}(\nabla_{\beta}\nabla_{\bar{\beta}}+\nabla_{\bar{\beta}}\nabla_{\beta})$.
}\\

Particularly, on Calabi-Yau manifolds, the flow can be reduced to a
hyperbolic Monge-Amp$\grave{e}$re equation:
\\$$\left\{\begin{array}{ll}
   \frac{\partial^2\varphi}{\partial t^2}(z,t)=\log\det g_{\alpha\bar{\beta}}(z,t)-\log\det g_{\alpha\bar{\beta}}(z,0)-f_0\\
   \varphi(z,0)=\varphi_0(z),\qquad\frac{\partial\varphi}{\partial
   t}(z,0)=\varphi_1(z)
  \end{array}\right.\eqno(3.2)$$\\
where
\\$$g_{\alpha\bar{\beta}}(z,t)\doteqdot g_{\alpha\bar{\beta}}(z,0)
+\frac{\partial^2\varphi}{\partial z^{\alpha}\partial \bar{z}^{\beta}}(z,t).$$\\

The paper is organized as follows. In Section 2, I will review
briefly some basics in K$\ddot{a}$hler geometry. In Section 3, I
will discuss the derivation and basic facts about hyperbolic
K$\ddot{a}$hler-Ricci flow in detail. In Section 4, evolutions of
various geometric quantities including the curvatures will be
derived.\\

\vspace{3mm}\par\noindent {\bf Acknowledgement. } {\it The author
thanks Professor De-Xing Kong and Professor Ke-Feng Liu for their
helpful suggestion without whose help, the come out of this paper is
impossible. The author also thanks Wen-Rong Dai and Xiao-Feng Sun
for their useful discussion.}
\vspace{3mm}\\

\section{Basic K${\bf\ddot{a}}$hler geometry}

I state some basic facts about K${\bf\ddot{a}}$hler geometry that
will be used in this paper. For a thorough study of
K${\bf\ddot{a}}$hler geometry, \cite{KN} \cite{MA} \cite{T}
\cite{ZF} are good choices. Let $\mathcal{M}^n$ be an n-dimensional
compact K$\ddot{a}$hler manifold. A Hermitian metric is given by
\\$$g=\sum g_{\alpha\bar{\beta}}dz^{\alpha}\otimes d\bar{z}^{\beta}.$$\\
Its associated K$\ddot{a}$hler form is
\\$$\omega=-\frac{1}{2}Im g=\frac{\sqrt{-1}}{2}\sum g_{\alpha\bar{\beta}}dz^{\alpha}\wedge d\bar{z}^{\beta}.$$\\
The K$\ddot{a}$hler condition requires that $\omega$ is a closed
positive global $(1,1)$-form.

The Riemann curvature is locally given by
\\$$R_{i\bar{j}k\bar{l}}=-\frac{\partial^2g_{i\bar{j}}}{\partial z^k\partial\bar{z}^l}
+\sum_{p,q=1}^ng^{p\bar{q}}\frac{\partial g_{i\bar{q}}}{\partial z^k}\frac{\partial g_{p\bar{j}}}{\partial\bar{z}^l},$$\\
and has the symmetries
\\$$R_{i\bar{j}k\bar{l}}=R_{k\bar{l}i\bar{j}}=-R_{\bar{j}ik\bar{l}}=-R_{i\bar{j}\bar{l}k}.$$\\

The Ricci curvature is defined by
\\$$R_{i\bar{j}}=g^{k\bar{l}}R_{i\bar{j}k\bar{l}},$$\\
and locally we have
\\$$R_{i\bar{j}}=-\frac{\partial^2\log\det(g_{k\bar{l}})}{\partial z^i\partial\bar{z}^j}.$$\\
The associated Ricci form is
\\$$\rho=\frac{\sqrt{-1}}{2}\sum R_{i\bar{j}}dz^i\wedge d\bar{z}^j.$$\\
It is a real closed $(1,1)$-form and represents the first Chern
class. The scalar curvature is
\\$$R=g^{i\bar{j}}R_{i\bar{j}}.$$\\

Given a K$\ddot{a}$hler metric, there are some other important
quantities. The nonzero Christoffel symbols are given by
\\$$\Gamma^k_{ij}=\sum_{l=1}^ng^{k\bar{l}}\frac{\partial g_{i\bar{l}}}{\partial z^j}\qquad
and\qquad \Gamma_{\bar{i}\bar{j}}^{\bar{k}}=\sum_{l=1}^ng^{\bar{k}l}\frac{\partial g_{l\bar{i}}}{\partial\bar{z}^j}.$$\\
The volume form is
\\$$d\mu=\frac{\omega^n}{n!}=(\sqrt{-1})^n\det(g_{i\bar{j}})dz^1\wedge d\bar{z}^1\wedge\cdot\cdot\cdot\wedge dz^n\wedge d\bar{z}^n.$$\\
The average scalar curvature is
\\$$r=\frac{\int_{\mathcal {M}^n}Rd\mu}{\int_{\mathcal {M}^n}d\mu}.$$\\

Finally I give an important lemma which is a consequence of the
Hodge decomposition theorem.\\

{\bf Lemma 2.1($\partial\bar{\partial}-$Lemma). }{\it Let
$\mathcal{M}^n$ be a compact K$\ddot{a}$hler manifold. If $a$ is an
exact real (1,1)-form, then there exists a real-valued function
$\psi$ on $\mathcal{M}^n$ such that
$\sqrt{-1}\partial\bar{\partial}\psi=a$. That is,
\\$$\frac{\partial^2}{\partial z^{\alpha}\partial\bar{z}^{\beta}}\psi=a_{\alpha\bar{\beta}},$$\\
where $a=\sqrt{-1}a_{\alpha\bar{\beta}}dz^{\alpha}\wedge
d\bar{z}^{\beta}$ and
$\overline{a_{\alpha\bar{\beta}}}=a_{\beta\bar{\alpha}}$.}\\

{\bf Proof. } This is a standard result in the theory of
K$\ddot{a}$hler manifolds, see \cite{ZF} for example.
\begin{flushright}
$\Box$
\end{flushright}

\section{Hyperbolic K$\ddot{a}$hler-Ricci flow}

In this section I will discuss the hyperbolic K$\ddot{a}$hler-Ricci
flow and its equivalent formulation as a single hyperbolic
Monge-Amp$\grave{e}$re equation. For the K$\ddot{a}$hler-Ricci flow
theory, please see \cite{CB}. I also suggest reader to read
Hamilton's papers for more details about Ricci flow \cite{H1}
\cite{H2}. The book \cite{CK} is also a good choice.

Given a compact complex manifold $\mathcal {M}^n$, consider the
hyperbolic K$\ddot{a}$hler-Ricci flow equation
\\$$\frac{\partial^2}{\partial t^2}g_{i\bar{j}}=-R_{i\bar{j}},$$\\
for a 1-parameter family of K$\ddot{a}$hler metrics, which is
obtained from the Ricci flow by dropping the factor of 2 with
respect to the convention. On a Calabi-Yau manifold, the first Chern
class vanishes, so $[\rho]=0$. Under this condition, from the
$\partial\bar{\partial}-Lemma$, the Ricci tensor has a potential,
i.e.
\\$$R_{\alpha\bar{\beta}}=\nabla_{\alpha}\nabla_{\bar{\beta}}f,$$\\
where $f(z,t)$ is a function defined on the manifold.

Assume
\\$$g_{\alpha\bar{\beta}}(z,t)=g_{\alpha\bar{\beta}}(z,0)+\nabla_{\alpha}\nabla_{\bar{\beta}}\varphi(z,t),$$\\
where $\varphi(z,t)$ is a function on the manifold, since we have
\\$$R_{\alpha\bar{\beta}}(t)=
-\nabla_{\alpha}\nabla_{\bar{\beta}}\log\det\biggl(g_{\gamma\bar{\delta}}(0)+
\frac{\partial^2\varphi(t)}{\partial z^{\gamma}\partial \bar{z}^{\delta}}\biggr).$$\\
Hence
\\ \begin{eqnarray*}
    \nabla_{\alpha}\nabla_{\bar{\beta}}\frac{\partial^2\varphi(z,t)}{\partial t^2}&=&\frac{\partial^2}{\partial
    t^2}g_{\alpha\bar{\beta}}(z,t)=-R_{\alpha\bar{\beta}}(z,t)=(R_{\alpha\bar{\beta}}(z,0)-R_{\alpha\bar{\beta}}(z,t))
    -R_{\alpha\bar{\beta}}(z,0)\\
    &=&\nabla_{\alpha}\nabla_{\bar{\beta}}\log\frac{\det(g_{\gamma\bar{\delta}}^0+\nabla_{\gamma}\nabla_{\bar{\delta}}\varphi)}
    {\det g_{\gamma\bar{\delta}}^0}-\nabla_{\alpha}\nabla_{\bar{\beta}}f_0.
    \end{eqnarray*}\\
So on Calabi-Yau manifolds, the hyperbolic K$\ddot{a}$hler-Ricci
flow equation is equivalent to the following hyperbolic (scalar)
complex Monge-Amp$\grave{e}$re equation due to the maximum principle
for compact manifolds:
\\$$\frac{\partial^2}{\partial t^2}\varphi=\log\frac{\det\biggl(g_{\gamma\bar{\delta}}^0
+\frac{\partial^2\varphi}{\partial z^{\gamma}\partial
\bar{z}^{\delta}}\biggr)}{\det g_{\gamma\bar{\delta}}^0}-f_0+c_1(t)$$\\
for some function of time $c_1(t)$ satisfying the compatibility
condition
\\$$\int_{\mathcal{M}^n}\biggl(\frac{\partial^2\varphi}{\partial t^2}-f_0\biggr)d\mu=\exp(c_1(t))Vol(\mathcal{M}^n),$$\\
here the Volume $Vol(\mathcal{M}^n)$ is under the metric
$g_{\alpha\bar{\beta}}(z,t)$. Further, by a time-dependent
translation of $\varphi(z,t)$, we can drop the factor $c_1(t)$.

Further we have
\\$$\frac{\partial}{\partial t}g_{\alpha\bar{\beta}}(z,0)=\nabla_{\alpha}\nabla_{\bar{\beta}}\frac{\partial}{\partial t}\varphi(z,0),$$\\
So we can let $[\frac{\partial}{\partial t}\omega(z,0)]=0$.

From the hyperbolicity and short time existence and uniqueness of
our hyperbolic Monge-Amp$\acute{e}$re equation, we have proved the
Theorem 1.1.

Next let us consider the normalized flow which is also studied in
\cite{KL} in the real case. Choose the normalization factor
$\varphi=\varphi(t)$ (Note: this has nothing to do with the
$\varphi$ above),
\\$$\tilde g_{i\bar{j}}=\varphi^2g_{i\bar{j}}$$\\
such that
\\$$\int_{\mathcal {M}^n}d\widetilde V=1,$$\\
and choose a new time parameter
\\$$\tilde t=\int_{\mathcal {M}^n}\varphi(t)dt.$$\\
Noting that for the normalized metric $\tilde g_{i\bar{j}}$, we have
\\$$\tilde R_{i\bar{j}}=R_{i\bar{j}},\qquad \tilde R=\frac{1}{\varphi^2}R,\qquad
\tilde r=\frac{1}{\varphi^2}r.$$\\
Thus
\\$$\frac{\partial\tilde g_{i\bar{j}}}{\partial\tilde t}=\varphi\frac{\partial g_{i\bar{j}}}{\partial t}
+2\frac{d\varphi}{dt}g_{i\bar{j}},$$
\\ \begin{eqnarray*}\frac{\partial^2\tilde g_{i\bar{j}}}{\partial\tilde
t^2}&=&\frac{\partial^2g_{i\bar{j}}}{\partial
t^2}+3\biggl(\frac{d}{dt}\log\varphi\biggr)\frac{\partial
g_{i\bar{j}}}{\partial
t}+2\biggl(\frac{d}{dt}\log\varphi\biggr)\biggl(\frac{d}{dt}\log\frac{d\varphi}{dt}\biggr)g_{i\bar{j}}\\
&=&-\tilde
R_{i\bar{j}}+3\frac{1}{\varphi}\biggl(\frac{d}{dt}\log\varphi\biggr)\frac{\partial\tilde
g_{i\bar{j}}}{\partial\tilde
t}+2\frac{1}{\varphi^2}\biggl(\frac{d}{dt}\log\varphi\biggr)\biggl(\frac{d}{dt}\log\frac{d\varphi}{dt}
-3\frac{d}{dt}\log\varphi\biggr)\tilde g_{i\bar{j}}\\ &=&-\tilde
R_{i\bar{j}}+a\frac{\partial\tilde g_{i\bar{j}}}{\partial\tilde
t}+b\tilde g_{i\bar{j}},
\end{eqnarray*}
\\
where $a$ and $b$ are certain functions of $t$.\\

Next let us consider the following hyperbolic system
\\$$\left\{\begin{array}{ll}
   \frac{\partial^2\varphi}{\partial t^2}(z,t)=\log\frac{\det g_{\alpha\bar{\beta}}(z,t)}
   {\det g_{\alpha\bar{\beta}}(z,0)}-f_0\\
   \varphi(z,0)=\varphi_0(z),\qquad\frac{\partial\varphi}{\partial
   t}(z,0)=\varphi_1(z)
  \end{array}\right.\eqno(3.1)$$\\
where
\\$$g_{\alpha\bar{\beta}}(z,t)\doteqdot g_{\alpha\bar{\beta}}(z,0)
+\frac{\partial^2\varphi}{\partial z^{\alpha}\partial \bar{z}^{\beta}}(z,t).$$\\

Let
\\$$v(z,t)\doteqdot-\frac{\partial\varphi}{\partial t}(z,t).$$\\
We get
\\$$\frac{\partial^2v}{\partial z^{\alpha}\partial\bar{z}^{\beta}}
=-\frac{\partial}{\partial
t}\biggl(\frac{\partial^2\varphi}{\partial
z^{\alpha}\partial\bar{z}^{\beta}}\biggr)
=-\frac{\partial}{\partial t}g_{\alpha\bar{\beta}}.$$\\
So
\\ \begin{eqnarray*}
   \frac{\partial^2}{\partial t^2}v &=&-\frac{\partial}{\partial t}\biggl(\frac{\partial^2\varphi}{\partial
   t^2}\biggr)\\
   &=&-\frac{\partial}{\partial t}(\log\det g_{\alpha\bar{\beta}})
   +\frac{\partial}{\partial t}(\log\det g_{\alpha\bar{\beta}}^0)-\frac{\partial}{\partial t}f_0\\
   &=&-g^{\alpha\bar{\beta}}\frac{\partial}{\partial
   t}g_{\alpha\bar{\beta}}\\
   &=&g^{\alpha\bar{\beta}}\frac{\partial^2 v}{\partial z^{\alpha}\partial
   \bar{z}^{\beta}}\\
   &=&\triangle v,
    \end{eqnarray*}\\
thus we have
\\$$\frac{\partial^2}{\partial t^2}v=\triangle v.$$\\
Note that the Laplacian operator here is time-dependent so the
equation is genuinely nonlinear.\\

\section{Evolutions of geometric quantities}

The hyperbolic K$\ddot{a}$hler-Ricci flow is a hyperbolic evolution
equation on the metrics. The evolution of the metrics implies
nonlinear wave equations for the Riemannian curvature tensor
$R_{i\bar{j}k\bar{l}}$, the Ricci curvature $R_{i\bar{j}}$ and the
scalar curvature $R$ which I will derive. I also derive the
evolutions of some other quantities in this section. For the
evolutions which are similar to those in this section associated to
K$\ddot{a}$hler-Ricci flow, please cite \cite{CB}.

Let $\mathcal{M}^n$ be an n-dimentional compact K$\ddot{a}$hler
manifold. Let us consider the hyperbolic K$\ddot{a}$hler-Ricci flow
on $\mathcal{M}^n$,
\\$$\frac{\partial^2}{\partial t^2}g_{i\bar{j}}(z,t)=-R_{i\bar{j}}(z,t).$$\\

{\bf Proposition 4.1(Riemannian curvature tensor). }{\it In a
unitary frame,
\\ \begin{eqnarray*}
  \frac{\partial^2}{\partial t^2}R_{i\bar{j}k\bar{l}}&=&\triangle_R R_{i\bar{j}k\bar{l}}
+R_{i\bar{\alpha}\beta\bar{l}}R_{\alpha\bar{j}k\bar{\beta}}-R_{i\bar{\alpha}k\bar{\beta}}R_{\alpha\bar{j}\beta\bar{l}}
+R_{i\bar{j}\beta\bar{\alpha}}R_{\alpha\bar{\beta}k\bar{l}}\\
&
&-\frac{1}{2}(R_{i\bar{\alpha}}R_{\alpha\bar{j}k\bar{l}}+R_{\alpha\bar{j}}R_{i\bar{\alpha}k\bar{l}}
+R_{k\bar{\alpha}}R_{i\bar{j}\alpha\bar{l}}+R_{\alpha\bar{l}}R_{i\bar{j}k\bar{\alpha}})\\
& &+2g^{p\bar{q}}\nabla_{k}\biggl(\frac{\partial}{\partial
t}g_{i\bar{q}}\biggr)\nabla_{\bar{l}}\biggl(\frac{\partial}{\partial
t}g_{p\bar{j}}\biggr)
    \end{eqnarray*}\\
The above formula also holds in arbitrary holomorphic coordinates if
repeated indices are contracted via the metric. I.e., if
$R_{i\bar{\alpha}\beta\bar{l}}R_{\alpha\bar{j}k\bar{\beta}}$ is
replaced by
$g^{\gamma\bar{\alpha}}g^{\beta\bar{\delta}}R_{i\bar{\alpha}\beta\bar{l}}R_{\gamma\bar{j}k\bar{\delta}}$,
etc.}\\

In the K$\ddot{a}$hler case, it is sometimes convenient to compute
locally in holomorphic coordinates in terms of ordinary derivatives.
The following lemma translates these ordinary derivatives to the
covariant derivatives, \cite{CB}.\\

{\bf Lemma 4.2(Relation between ordinary and covariant derivatives).
}{\it Let $\eta$ be a closed (1,1)-form. Locally it is represented
by $\eta_{\alpha\bar{\beta}}$, which is Hermitian symmetric. Let
$\eta_{\alpha\bar{\beta},\gamma\bar{\delta}}$ denote the covariant
derivatives and $\eta_{\alpha\bar{\beta}\gamma\bar{\delta}}$ denote
$\frac{\partial^2}{\partial z^{\gamma}\partial
z^{\bar{\delta}}}\eta_{\alpha\bar{\beta}}$. Then at the center $x$
of normal holomorphic coordinates
\\$$\eta_{\gamma\bar{\delta},\alpha\bar{\beta}}=\eta_{\gamma\bar{\delta}\alpha\bar{\beta}}
+\eta_{s\bar{\delta}}R_{\alpha\bar{\beta}\gamma\bar{s}},$$
$$\eta_{\gamma\bar{\delta},\bar{\beta}\alpha}=\eta_{\gamma\bar{\delta}\bar{\beta}\alpha}
+\eta_{\gamma\bar{s}}R_{\alpha\bar{\beta}s\bar{\delta}}.$$\\}

{\bf Proof. } \begin{eqnarray*}
\eta_{\gamma\bar{\delta},\alpha\bar{\beta}}&=&\nabla_{\bar{\beta}}\nabla_{\alpha}\eta_{\gamma\bar{\delta}}
=\partial_{\bar{\beta}}\nabla_{\alpha}\eta_{\gamma\bar{\delta}}-
\bar{\Gamma}^{\varepsilon}_{\beta\delta}\nabla_{\alpha}\eta_{\gamma\bar{\varepsilon}}\\
&=&\partial_{\bar{\beta}}(\partial_{\alpha}\eta_{\gamma\bar{\delta}}-
\Gamma^{\varepsilon}_{\alpha\gamma}\eta_{\varepsilon\bar{\delta}})\\
&=&\partial_{\bar{\beta}}\partial_{\alpha}\eta_{\gamma\bar{\delta}}-
(\partial_{\bar{\beta}}\Gamma^{\varepsilon}_{\alpha\gamma})\eta_{\varepsilon\bar{\delta}}
-\Gamma^{\varepsilon}_{\alpha\gamma}(\partial_{\bar{\beta}}\eta_{\varepsilon\bar{\delta}})\\
&=&\eta_{\gamma\bar{\delta}\alpha\bar{\beta}}+R^{\varepsilon}_{\alpha\bar{\beta}\gamma}\eta_{\varepsilon\bar{\delta}}.
\end{eqnarray*}\\
In the above calculation, we use the fact that $\eta$ is closed so
$\nabla_{\alpha}\eta_{\gamma\bar{\varepsilon}}=\partial_{\bar{\beta}}\eta_{\varepsilon\bar{\delta}}=0$.
The second formula is the conjugate of the first.

\begin{flushright}
$\Box$
\end{flushright}

{\bf Proof of Proposition 4.1. }Recall that
\\$$R_{i\bar{j}k\bar{l}}=-\frac{\partial^2g_{i\bar{j}}}{\partial z^k\partial \bar{z}^l}
+\sum_{p,q=1}^ng^{p\bar{q}}\frac{\partial g_{i\bar{q}}}{\partial z^k}\frac{\partial g_{p\bar{j}}}{\partial \bar{z}^l}$$\\
This implies that, in a normal holomorphic coordinate system
centered at any given point
\\ \begin{eqnarray*}
   \frac{\partial^2}{\partial t^2}R_{i\bar{j}k\bar{l}} &=&
   -\frac{\partial^2}{\partial z^k\partial\bar{z}^l}\biggl(\frac{\partial^2}{\partial
   t^2}g_{i\bar{j}}\biggr)+2g^{p\bar{q}}\nabla_{k}\biggl(\frac{\partial}{\partial
t}g_{i\bar{q}}\biggr)\nabla_{\bar{l}}\biggl(\frac{\partial}{\partial
t}g_{p\bar{j}}\biggr)\\
   &=&\frac{\partial^2}{\partial
   z^k\partial\bar{z}^l}R_{i\bar{j}}+2g^{p\bar{q}}\nabla_{k}\biggl(\frac{\partial}{\partial
t}g_{i\bar{q}}\biggr)\nabla_{\bar{l}}\biggl(\frac{\partial}{\partial
t}g_{p\bar{j}}\biggr)\\
   &=&\nabla_k\nabla_{\bar{l}}R_{i\bar{j}}-R_{i\bar{\alpha}}R_{\alpha\bar{j}k\bar{l}}
   +2g^{p\bar{q}}\nabla_{k}\biggl(\frac{\partial}{\partial
t}g_{i\bar{q}}\biggr)\nabla_{\bar{l}}\biggl(\frac{\partial}{\partial
t}g_{p\bar{j}}\biggr).
    \end{eqnarray*}\\
The proposition follows from the following two identities:
\\ \begin{eqnarray*}
   \nabla_k\nabla_{\bar{l}}R_{i\bar{j}} &=& \nabla_k\nabla_{\bar{l}}R_{i\bar{j}\beta\bar{\beta}}
   = \nabla_k\nabla_{\bar{\beta}}R_{i\bar{j}\beta\bar{l}}\\
   &=&\nabla_{\bar{\beta}}\nabla_{\beta}R_{i\bar{j}k\bar{l}}
   -R_{i\bar{\alpha}k\bar{\beta}}R_{\alpha\bar{j}\beta\bar{l}}
   +R_{\alpha\bar{j}k\bar{\beta}}R_{i\bar{\alpha}\beta\bar{l}}\\
   & &-R_{\beta\bar{\alpha}k\bar{\beta}}R_{i\bar{j}\alpha\bar{l}}+R_{\alpha\bar{l}k\bar{\beta}}R_{i\bar{j}\beta\bar{\alpha}},
    \end{eqnarray*}\\
and
\\ \begin{eqnarray*}
    \nabla_{\bar{\beta}}\nabla_{\beta}R_{i\bar{j}k\bar{l}} &=& \nabla_{\beta}\nabla_{\bar{\beta}}R_{i\bar{j}k\bar{l}}+
    R_{i\bar{\alpha}\beta\bar{\beta}}R_{\alpha\bar{j}k\bar{l}}
    -R_{\alpha\bar{j}\beta\bar{\beta}}R_{i\bar{\alpha}k\bar{l}}\\
    & &+R_{k\bar{\alpha}\beta\bar{\beta}}R_{i\bar{j}\alpha\bar{l}}
    -R_{\alpha\bar{l}\beta\bar{\beta}}R_{i\bar{j}k\bar{\alpha}}.
    \end{eqnarray*}\\
Recall that
$\triangle_R=\frac{1}{2}(\nabla_{\beta}\nabla_{\bar{\beta}}+\nabla_{\bar{\beta}}\nabla_{\beta})$.
\begin{flushright}
$\Box$
\end{flushright}

{\bf Lemma 4.3. }{\it We will need the following two basic formulas
under hyperbolic K$\ddot{a}$hler-Ricci flow:
\\$$\frac{\partial}{\partial t}g^{k\bar{l}}=-g^{k\bar{s}}\biggl(\frac{\partial}{\partial t}g_{r\bar{s}}\biggr)g^{r\bar{l}},$$\\
\\$$\frac{\partial^2}{\partial
t^2}g^{k\bar{l}}=R_{r\bar{s}}g^{k\bar{s}}g^{r\bar{l}}
+2\biggl(\frac{\partial}{\partial
t}g_{n\bar{m}}\biggr)\biggl(\frac{\partial}{\partial
t}g_{r\bar{s}}\biggr) g^{k\bar{s}}g^{r\bar{m}}g^{n\bar{l}}.$$\\}

{\bf Corollary 4.4(Ricci curvature). }{\it The Ricci curvature
satisfies the following in a unitary frame:
\\ \begin{eqnarray*}\frac{\partial^2}{\partial
t^2}R_{i\bar{j}}&=&\triangle_R
R_{i\bar{j}}+R_{i\bar{j}k\bar{l}}R_{l\bar{k}}
-R_{i\bar{k}}R_{k\bar{j}}-2\biggl\langle\frac{\partial}{\partial
t}g_{k\bar{l}},\frac{\partial}{\partial
t}R_{i\bar{j}k\bar{l}}\biggr\rangle\\
& &+2R_{i\bar{j}k\bar{l}}\biggl(\frac{\partial}{\partial
t}g_{n\bar{m}}\biggr)\biggl(\frac{\partial}{\partial
t}g_{r\bar{s}}\biggr)g^{r\bar{m}}g^{n\bar{l}}g^{k\bar{s}}\\
&
&+2g^{p\bar{q}}g^{k\bar{l}}\nabla_{k}\biggl(\frac{\partial}{\partial
t}g_{i\bar{q}}\biggr)\nabla_{\bar{l}}\biggl(\frac{\partial}{\partial
t}g_{p\bar{j}}\biggr)\end{eqnarray*}}

{\bf Proof. }
\\ \begin{eqnarray*}
    \frac{\partial^2}{\partial t^2}R_{i\bar{j}}&=&\frac{\partial^2}{\partial
    t^2}(g^{k\bar{l}}R_{i\bar{j}k\bar{l}})\\
    &=&g^{k\bar{l}}\biggl(\frac{\partial^2}{\partial
    t^2}R_{i\bar{j}k\bar{l}}\biggr)+R_{i\bar{j}k\bar{l}}\biggl(\frac{\partial^2}{\partial
    t^2}g^{k\bar{l}}\biggr)+2\biggl(\frac{\partial}{\partial t}g^{k\bar{l}}\biggr)
    \biggl(\frac{\partial}{\partial t}R_{i\bar{j}k\bar{l}}\biggr).
    \end{eqnarray*}\\
Put the evolutions of $R_{i\bar{j}k\bar{l}}$ and $g^{k\bar{l}}$ in,
we get the result.
\begin{flushright}
$\Box$
\end{flushright}

{\bf Proposition 4.5(Evolution of R). }{\it The scalar curvature $R$
evolves by
\\ \begin{eqnarray*}\frac{\partial^2}{\partial t^2}R&=&\triangle R+|R_{\alpha\bar{\beta}}|^2
+\triangle\biggl|\frac{\partial}{\partial
t}g_{\alpha\bar{\beta}}\biggr|^2
-2\biggl\langle\frac{\partial}{\partial
t}g_{\alpha\bar{\beta}},\frac{\partial}{\partial
t}R_{\alpha\bar{\beta}}\biggr\rangle\\
& &+2R_{k\bar{l}}\biggl(\frac{\partial}{\partial
t}g_{n\bar{m}}\biggr)\biggl(\frac{\partial}{\partial
t}g_{r\bar{s}}\biggr)g^{r\bar{m}}g^{n\bar{l}}g^{k\bar{s}}.\end{eqnarray*}\\}

{\bf Proof. }First we have
\\$$\frac{\partial}{\partial t}\log\det g=g^{\alpha\bar{\beta}}\frac{\partial}{\partial t}g_{\alpha\bar{\beta}},$$\\
\\ \begin{eqnarray*}
   \frac{\partial^2}{\partial t^2}\log\det g &=&-\biggl|\frac{\partial}{\partial t}g\biggr|^2
   +g^{\alpha\bar{\beta}}\frac{\partial^2}{\partial
   t^2}g_{\alpha\bar{\beta}}\\
   &=&-\biggl|\frac{\partial}{\partial t}g\biggr|^2-R.
    \end{eqnarray*}\\
Thus
\\$$\frac{\partial^2}{\partial t^2}R_{\alpha\bar{\beta}}=
-\nabla_{\alpha}\nabla_{\bar{\beta}}\biggl(\frac{\partial^2}{\partial
t^2}\log\det g\biggr)=
\nabla_{\alpha}\nabla_{\bar{\beta}}R+\nabla_{\alpha}\nabla_{\bar{\beta}}\biggl|\frac{\partial}{\partial t}g\biggr|^2.$$\\
So finally
\\ \begin{eqnarray*}
    \frac{\partial^2}{\partial t^2}R &=&
     -2\biggl\langle\frac{\partial}{\partial t}g_{\gamma\bar{\delta}},\frac{\partial}{\partial t}
     R_{\gamma\bar{\delta}}\biggr\rangle
     +g^{\alpha\bar{\beta}}\frac{\partial^2}{\partial t^2}R_{\alpha\bar{\beta}}
     +\biggl(\frac{\partial^2}{\partial t^2}g^{\alpha\bar{\beta}}\biggr)R_{\alpha\bar{\beta}}\\
     &=&\triangle R+|R_{\alpha\bar{\beta}}|^2
+\triangle\biggl|\frac{\partial}{\partial
t}g_{\alpha\bar{\beta}}\biggr|^2
-2\biggl\langle\frac{\partial}{\partial
t}g_{\alpha\bar{\beta}},\frac{\partial}{\partial
t}R_{\alpha\bar{\beta}}\biggr\rangle\\
& &+2R_{k\bar{l}}\biggl(\frac{\partial}{\partial
t}g_{n\bar{m}}\biggr)\biggl(\frac{\partial}{\partial
t}g_{r\bar{s}}\biggr)g^{r\bar{m}}g^{n\bar{l}}g^{k\bar{s}}
    \end{eqnarray*}\\
\begin{flushright}
$\Box$
\end{flushright}

Next I present the evolution equations of the Christoffel symbols
and the volume form.\\

{\bf Proposition 4.6. }{\it The christoffel symbols and the volume
form evolve respectively by
\\$$\frac{\partial^2}{\partial t^2}\Gamma_{\alpha\beta}^{\gamma}=
-g^{\gamma\bar{\delta}}\nabla_{\alpha}R_{\beta\bar{\delta}}+
2\biggl(\frac{\partial}{\partial
t}g^{\gamma\bar{\delta}}\biggr)\nabla_{\alpha}\biggl(\frac{\partial}
{\partial t}g_{\beta\bar{\delta}}\biggr),$$\\
\\$$\frac{\partial^2}{\partial t^2}d\mu=\biggl[-R+\biggl(g^{\alpha\bar{\beta}}
\frac{\partial g_{\alpha\bar{\beta}}}{\partial t}\biggr)^2
-\biggl|\frac{\partial g_{\alpha\bar{\beta}}}{\partial
t}\biggr|^2\biggr]d\mu.$$\\}

{\bf Proof. }Recall that
\\$$\Gamma_{\alpha\beta}^{\gamma}=\frac{1}{2}g^{\gamma\bar{\delta}}\biggl(\frac{\partial}{\partial z^{\alpha}}g_{\beta\bar{\delta}}
+\frac{\partial}{\partial
z^{\beta}}g_{\alpha\bar{\delta}}-\frac{\partial}{\partial\bar{z}^{\delta}}g_{\alpha\beta}\biggr)
=g^{\gamma\bar{\delta}}\frac{\partial}{\partial z^{\alpha}}g_{\beta\bar{\delta}},$$\\
thus
\\$$\frac{\partial}{\partial t}\Gamma_{\alpha\beta}^{\gamma}=g^{\gamma\bar{\delta}}\frac{\partial}{\partial z^{\alpha}}
\biggl(\frac{\partial}{\partial t}g_{\beta\bar{\delta}}\biggr)
+\biggl(\frac{\partial}{\partial t}g^{\gamma\bar{\delta}}\biggr)\frac{\partial}{\partial z^{\alpha}}g_{\beta\bar{\delta}}.$$\\
So finally in normal coordinate
\\ \begin{eqnarray*}
    \frac{\partial^2}{\partial t^2}\Gamma_{\alpha\beta}^{\gamma} &=&
    2\biggl(\frac{\partial}{\partial t}g^{\gamma\bar{\delta}}\biggr)\frac{\partial}
    {\partial z^{\alpha}}\biggl(\frac{\partial}{\partial t}g_{\beta\bar{\delta}}\biggr)
    +g^{\gamma\bar{\delta}}\frac{\partial}{\partial z^{\alpha}}\biggl(\frac{\partial^2}{\partial
    t^2}g_{\beta\bar{\delta}}\biggr)+\biggl(\frac{\partial^2}
    {\partial t^2}g^{\gamma\bar{\delta}}\biggr)\biggl(\frac{\partial}{\partial z^{\alpha}}g_{\beta\bar{\delta}}\biggr)\\
   &=&-g^{\gamma\bar{\delta}}\nabla_{\alpha}R_{\beta\bar{\delta}}+
    2\biggl(\frac{\partial}{\partial t}g^{\gamma\bar{\delta}}\biggr)\nabla_{\alpha}
    \biggl(\frac{\partial}{\partial t}g_{\beta\bar{\delta}}\biggr).
    \end{eqnarray*}\\

The evolution of $d\mu$ follows from the quantities
\\$$\frac{\partial}{\partial t}\log\det g=g^{\alpha\bar{\beta}}\frac{\partial}{\partial t}g_{\alpha\bar{\beta}}$$\\
and
\\$$\frac{\partial^2}{\partial t^2}\det g_{\alpha\bar{\beta}}=
\frac{\partial}{\partial t}(\det g)\cdot
g^{\alpha\bar{\beta}}\frac{\partial}{\partial
t}g_{\alpha\bar{\beta}} -(\det g)\biggl|\frac{\partial}{\partial
t}g_{\alpha\bar{\beta}}\biggr|^2+(\det g)\cdot
g^{\alpha\bar{\beta}}(-R_{\alpha\bar{\beta}}).$$
\begin{flushright}
$\Box$
\end{flushright}

{\bf Proposition 4.7(potential of Ricci curvature). }{\it On compact
K$\ddot{a}$hler manifolds, the potential $f$ of Ricci curvature
$(R_{\alpha\bar{\beta}}=\nabla_{\alpha}\nabla_{\bar{\beta}}f)$
satisfies:
\\$$\frac{\partial^2}{\partial t^2}f=\triangle f+\biggl|\frac{\partial}{\partial
t}g\biggr|^2+c(t).$$\\}

{\bf Proof. }
\\ \begin{eqnarray*}
\nabla_{\alpha}\nabla_{\bar{\beta}}\biggl(\frac{\partial^2}{\partial
t^2}f\biggr)&=&\frac{\partial^2}{\partial
t^2}\nabla_{\alpha}\nabla_{\bar{\beta}}f=\frac{\partial^2}{\partial
t^2}R_{\alpha\bar{\beta}}\\
&=&\nabla_{\alpha}\nabla_{\bar{\beta}}\biggl(R+\biggl|\frac{\partial}{\partial
t}g\biggr|^2\biggr)\\
&=&\nabla_{\alpha}\nabla_{\bar{\beta}}\biggl(\triangle
f+\biggl|\frac{\partial}{\partial t}g\biggr|^2\biggr),
\end{eqnarray*}\\
so from the maximum principle of compact manifolds,
\\$$\frac{\partial^2}{\partial t^2}f=\triangle f+\biggl|\frac{\partial}{\partial
t}g\biggr|^2+c(t).$$\\

\begin{flushright}
$\Box$
\end{flushright}

By a suitable time-dependent translation of $f$, we can drop the
factor $c(t)$.

{\bf Proposition 4.8. }{\it On compact K$\ddot{a}$hler manifolds,
the evolutions of $\int_{\mathcal {M}^n}Rd\mu$ and $r$ are:
\\ \begin{eqnarray*}\frac{d^2}{dt^2}\int_{\mathcal{M}^n}Rd\mu&=&\int_{\mathcal {M}^n}R\biggl[\biggl(g^{\alpha\bar{\beta}}
\frac{\partial g_{\alpha\bar{\beta}}}{\partial
t}\biggr)^2-\biggl|\frac{\partial g_{\alpha\bar{\beta}}}{\partial
t}\biggr|^2\biggr]d\mu +2\int_{\mathcal
{M}^n}\frac{\partial}{\partial t}R\frac{\partial}{\partial
t}d\mu\\
& &-2\int_{\mathcal {M}^n}\biggl\langle\frac{\partial}{\partial t}g_{\alpha\bar{\beta}},
\frac{\partial}{\partial t}R_{\alpha\bar{\beta}}\biggr\rangle d\mu\\
& &+2\int_{\mathcal
{M}^n}R_{k\bar{l}}\biggl(\frac{\partial}{\partial
t}g_{n\bar{m}}\biggr)\biggl(\frac{\partial}{\partial
t}g_{r\bar{s}}\biggr)g^{r\bar{m}}g^{n\bar{l}}g^{k\bar{s}} d\mu,\end{eqnarray*}\\

\begin{eqnarray*}
    \frac{d^2r}{dt^2}&=&r^2+\frac{\biggl(\frac{d^2}{dt^2}\int_{\mathcal {M}^n}Rd\mu\biggr)}{\biggl(\int_{\mathcal
    {M}^n}d\mu\biggr)}+\frac{2\biggl(\frac{d}{dt}\int_{\mathcal {M}^n}d\mu\biggr)^2\biggl(\int_{\mathcal {M}^n}Rd\mu\biggr)}
    {\biggl(\int_{\mathcal {M}^n}d\mu\biggr)^3}\\
    & &-\frac{\biggl(\int_{\mathcal {M}^n}Rd\mu\biggr)\biggl(\int_{\mathcal {M}^n}\biggl(g^{\alpha\bar{\beta}}\frac{\partial g_{\alpha\bar{\beta}}}{\partial t}\biggr)^2-
    \biggl|\frac{\partial g_{\alpha\bar{\beta}}}{\partial t}\biggr|^2d\mu\biggr)
    +2\biggl(\frac{d}{dt}\int_{\mathcal {M}^n}Rd\mu\biggr)\biggl(\frac{d}{dt}\int_{\mathcal {M}^n}d\mu\biggr)}{\biggl(\int_{\mathcal
    {M}^n}d\mu\biggr)^2}.
    \end{eqnarray*}\\
}

{\bf Proof. }Since
\\$$\frac{d^2}{dt^2}\int_{\mathcal{M}^n}Rd\mu=\int_{\mathcal
{M}^n}\frac{\partial^2}{\partial t^2}Rd\mu+\int_{\mathcal
{M}^n}R\frac{\partial^2}{\partial t^2}d\mu
+2\int_{\mathcal {M}^n}\frac{\partial}{\partial t}R\frac{\partial}{\partial t}d\mu,$$\\
put the evolutions of $R$ and $\mu$ in, we get the result for
$\int_{\mathcal {M}^n}Rd\mu$.\\

Since
\begin{eqnarray*} \frac{dr}{dt}=\frac{\biggl(\frac{d}{dt}\int_{\mathcal
{M}^n}R d\mu\biggr)\biggl(\int_{\mathcal {M}^n}d\mu\biggr)
-\biggl(\int_{\mathcal {M}^n}R
d\mu\biggr)\biggl(\frac{d}{dt}\int_{\mathcal {M}^n}d\mu\biggr)}
{\biggl(\int_{\mathcal {M}^n}d\mu\biggr)^2},\end{eqnarray*}\\

and
\begin{eqnarray*}\frac{d^2r}{dt^2}&=&\frac{\biggl(\frac{d^2}{dt^2}\int_{\mathcal {M}^n}Rd\mu\biggr)}
{\biggl(\int_{\mathcal
{M}^n}d\mu\biggr)}-\frac{\biggl(\frac{d^2}{dt^2}\int_{\mathcal
{M}^n}d\mu\biggr)\biggl(\int_{\mathcal
{M}^n}Rd\mu\biggr)+2\biggl(\frac{d}{dt}\int_{\mathcal
{M}^n}d\mu\biggr)\biggl(\frac{d}{dt}\int_{\mathcal
{M}^n}Rd\mu\biggr)}{\biggl(\int_{\mathcal {M}^n}d\mu\biggr)^2}\\
&&+\frac{2\biggl(\frac{d}{dt}\int_{\mathcal
{M}^n}d\mu\biggr)^2\biggl(\int_{\mathcal
{M}^n}Rd\mu\biggr)}{\biggl(\int_{\mathcal {M}^n}d\mu\biggr)^3},
\end{eqnarray*}\\

put the evolutions of $R$ and $\mu$ in, we get the result for $r$.

\begin{flushright}
$\Box$
\end{flushright}

{\bf Remark. } In K$\ddot{a}$hler-Ricci flow, the evolutions of
$\int_{\mathcal {M}^n}Rd\mu$ and $r$ are
\\$$\frac{\partial}{\partial t}\int_{\mathcal {M}^n}Rd\mu=0,$$\\
and
\\$$\frac{\partial}{\partial t}r=r^2.$$\\
While in our case, all the extra terms are of first-order
derivatives.

\section{Further discussions}

The flow considered here is the complex version of the hyperbolic
geometric flow introduced by Kong and Liu \cite{KL}. Kong et al
\cite{KLX} proved that on Riemann surfaces, the long time existence
of hyperbolic geometric flow depends on the choice of the initial
velocity. We can expect that its complex version has similar
property. On Calabi-Yau manifolds, the flow can be simplified to a
single complex hyperbolic Monge-Amp$\grave{e}$re equation. Note that
Yau \cite{Y} used the elliptic Monge-Amp$\grave{e}$re equation to
prove the famous Calabi's conjecture and Cao \cite{C} used its
parabolic version and techniques from K$\ddot{a}$hler-Ricci flow to
reprove Calabi's conjecture. I hope the hyperbolic
Monge-Amp$\grave{e}$re equation is also powerful to understand
K$\ddot{a}$hler manifolds. We can expect this new flow is helpful to
study wave phenomena in the nature especially the Einstein equation.
In the future, we will study several fundamental problems on the
hyperbolic K$\ddot{a}$hler-Ricci flow, for example, long-time
existence, formation of singularities as well as physical
applications.

\vskip 10mm
\noindent

\end{document}